\documentclass[a4paper,12pt]{article}
\usepackage[latin1]{inputenc}
\usepackage{amssymb}
\usepackage{amsmath}
\usepackage{color}

\usepackage{cancel}

\newtheorem{theorem}{Theorem}

\newtheorem{definition}[theorem]{Definition}

\newtheorem{corollary}[theorem]{Corollary}
\newtheorem{example}[theorem]{Example}
\newtheorem{lemma}[theorem]{Lemma}

\newtheorem{remark}[theorem]{Remark}

\def\trans#1                            %
    {{                          %
    \includegraphics{/home/harald/texinput/harald/trans.ps}    %
    \phantom{tt}                        %
    }}                          %
\def\transversal{{\pitchfork\kern -.907 em \overline{\phantom{\mbox{\footnotesize\rm lr}}}\kern .35em}}

\def\BA{\begin{array}} \def\EA{\end{array}}             
\def\BI{\begin{itemize}} \def\EI{\end{itemize}}             
\def\BM{\left(\begin{array}} \def\EM{\end{array}\right)}        
\def\BSet{\left\{\begin{array}{c|c}} \def\ESet{\end{array}\right\}} 

\def\ds{\displaystyle}

\def\QED{{\hfill $\square$}}
\def\pphi{{\varphi}}
\def\veps{\varepsilon}

\def\N{{\mathbb N}}

\def\R{{\mathbb R}}

\def\proof{{\noindent{\em Proof. }}}

\def\argmin{\mathop{{\rm argmin}}}

\def\ind{{\rm i\kern -.4pt n\kern -.4pt d}}

\def\span{\mathop{{\rm span}}}

\def\cala{{\cal A}}
\def\calb{{\cal B}}
\def\calc{{\cal C}}

\def\calf{{\cal F}}

\def\cali{{\cal I}}

\def\calp{{\cal P}}

\def\?{\textquestiondown}
\def\!{\textexclamdown}

\newcommand{\SP}{{\mathop{SP}}}
\newcommand{\MF}{{\mathop{MF}}}
\newcommand{\SQP}{{\mathop{SQP}}}
\newcommand{\jsp}{{\mathop{jSP}}}
\newcommand{\JSP}{{\mathop{JSP}}}
\newcommand{\jmf}{{\mathop{jMF}}}
\newcommand{\JMF}{{\mathop{JMF}}}
\newcommand{\graph}{{\mathop{graph}}}
\newcommand{\sigmasp}{{\sigma\mathop{SP}}}
\newcommand{\sigmamf}{{\sigma\mathop{MF}}}
\newcommand{\SLACK}{{\mathop{SLACK}}}
\newcommand{\SPMF}{{\mathop{SP2MF}}}
\newcommand{\SPMFtar}{{\mathop{SP2MF}^{\mathop{tar}}}}
\newcommand{\SPMFvar}{{\mathop{SP2MF}^{\mathop{var}}}}
\newcommand{\SPMFcode}{{\mathop{SP2MF}^{\mathop{code}}}}

\newcommand{\SPMFfunc}{{\mathop{SP2MF}^{\mathop{func}}}}
\newcommand{\MFSP}{{\mathop{MF2SP}}}
\newcommand{\MFSPtar}{{\mathop{MF2SP}^{\mathop{tar}}}}
\newcommand{\MFSPvar}{{\mathop{MF2SP}^{\mathop{var}}}}
\newcommand{\MFSPfunc}{{\mathop{MF2SP}^{\mathop{func}}}}

\newcommand{\chii}{{\raisebox{.5ex}{$\chi$}}}
\newcommand{\id}{{\mathop{id}}}

\newcommand{\aeq}[2]{{\kern .1em {}_{#1} \kern -.4em \sim \kern -.4em {}_{#2} \kern .1em}}
\newcommand{\simsim}{{\kern .1em \sim \kern .1em}}

\newcommand{\Ttar}{{T^{\mathop{tar}}}}
\newcommand{\Tset}{{T^{\mathop{set}}}}
\newcommand{\Tvar}{{T^{\mathop{var}}}}
\newcommand{\Tfunc}{{T^{\mathop{func}}}}
\newcommand{\Tcode}{{T^{\mathop{code}}}}
\newcommand{\TTtar}{{T_2^{\mathop{tar}}}}
\newcommand{\TTcode}{{T_2^{\mathop{code}}}}
\newcommand{\TTfunc}{{T_2^{\mathop{func}}}}

\newcommand{\TTTcode}{{T_3^{\mathop{code}}}}
\newcommand{\TTTfunc}{{T_3^{\mathop{func}}}}
\newcommand{\TTTTfunc}{{(T_3\circ T_2)^{\mathop{func}}}}
\newcommand{\TTTTcode}{{(T_3\circ T_2)^{\mathop{code}}}}

\newcommand{\ttar}{{T_1^{\mathop{tar}}}}
\newcommand{\tvar}{{T_1^{\mathop{tar}}}}

\newcommand{\tfunc}{{T_1^{\mathop{func}}}}

\title{Stationary Point Sets:\\
Convex Quadratic Optimization is Universal\\
in Nonlinear Optimization}
\author{Harald G\"unzel\thanks{
Department of Mathematics -- C,
RWTH Aachen University,
D-52056 Aachen,
Germany,
email: \texttt{guenzel@mathc.rwth-aachen.de}
}
}
\date{October 24, 2012, and August 05, 2013}
\begin{document}
\maketitle

\null\vfill
\begin{center}
Dedicated to Bert Jongen
\vspace{12pt}
\end{center}
\null\vfill

\begin{abstract}
We investigate the local topological structure,
stationary point sets in parametric optimization genericly may have.
Our main result states that, up to stratified isomorphism,
any such structure is already present in the small subclass
of parametric problems with convex quadratic objective function
and affine-linear constraints.

In other words,
the convex quadratic problems produce
a normal form for the local topological structure of stationary point sets.

As a consequence we see, as far as no equality constraints are involved,
that the closure of the stationary point set constitutes a manifold with boundary.
The boundary is exactly the violation set of the Mangasarian Fromovitz
constraint qualification.

A side result states that stationary point sets and violation sets
of Mangasarian Fromovitz constraint qualification
carry the same set of possible local structures as stratified spaces.
\end{abstract}

\section{Introduction}
We consider parametric families of constraint sets $M(y)$ and of optimization problems
$\calp(y)$. Here, $y\in\R^p$ is a parameter and the objects of our study are
defined as follows:
\begin{eqnarray}
\nonumber  M(y) &= \{x\in\R^n~|~&g_i(x,y)\le 0,\,i=1\dots m_\le,\\
                       &&h_j(x,y)=0,\,j=1\dots m_=\},\label{eq:M}\\
  \calp(y) &\mathop{minimize}~&f(x,y)~\mbox{ s.t. }~x\in M(y)\label{eq:calp},
\end{eqnarray}
where the appearing objective function $f$ and the constraint functions $g_i$, $h_j$
are supposed to be smooth, i.e. from $C^\infty(\R^{n+p})$.
Throughout the paper, for sake of a simpler exposition,
all mappings, manifolds, etc., are supposed to be smooth unless stated otherwise.
For various reasons it is convenient to refer to the objective function $f$ also
with the symbol $g_{m_\le^*}$, where $m_\le^*:=m_\le+1$,
i.e. we use the identity $g_{m_\le^*}\equiv f$.

We are mainly interested in the local topological structure of the stationary point set,
defined by setting
\begin{equation}\label{eq:SP}
    \SP:=\SP(f,g,h):=\{(x,y)~|~x\mbox{ is a stationary point for }\calp(y)\}.
\end{equation}
Knowledge of this structure may be helpful in various instances,
especially if $\SP$ serves as a constraint set of a so-called upper level problem,
such as in bilevel optimization or in mathematical programming with equilibrium constraints,
c.f. \cite{lit:equi}.

Our main result is some kind of universality theorem,
roughly stating that any topological structure of a stationary point set
appearing in general nonlinear optimization may already be found
in convex quadratic optimization.
This, for instance, offers the opportunity to compare the
topological structure of a stationary point set at unstable situations,
i.e. in the neighborhood of pairs $(x,y)$,
where $\SP$ cannot be locally written as the graph of a continuous function
$y\mapsto x(y)$, with stable ones, that are immanent in convex quadratic optimization.

In order to formulate our main result more precisely, we need some concepts.
First, we introduce a combinatorial code for feasible points $x\in M(y)$.
One (but not the main) reason
for using this concept
is the precise definition of stationary points
and the Mangasarian Fromovitz constraint qualification (MFCQ).
We start with the indexed Lagrange function, defined for $I\subset\{1\dots m_\le^*\}$,
$J\subset\{1\dots m_=\}$,
and operating on the variables
$(x,y)\in\R^{n+p}$, and $\mu_i,\lambda_j\in\R$
by setting
\begin{equation}\label{eq:L}
    L_{I,J}:=\sum_{i\in I}\mu_i g_i(x,y)+ \sum_{j\in J}\lambda_j h_j(x,y).
\end{equation}
The combinatorial code of a feasible point $x\in M(y)$ consists of two parts.
The first part is the index set $I_0=I_0(x,y)$ of {\em active} inequality constraints,
defined by setting
$I_0:=\{i\in\{1...m\}~|~g_i(x,y)=0\}$.
Moreover let $I^\ast_0:=I_0\cup\{m_\le^*\}$,
only defined in case that an objective function $f$ is present.
The second part of the code is the set $\cali=\cali(x,y)$,
containing all pairs of index sets $(I,J)$, $I\subset I_0^*$, $J\subset\{1\dots m_=\}$
whose disjoint union is minimal (w.r.t. inclusion)
within the set of all such pairs with the properties that their disjoint
union is nonempty and that there exist multipliers $\mu_i>0$, $i\in I$,
$\lambda_j\neq 0$, $j\in J$, such that $D_xL_{I,J}=0$.
Letting $\cali^\SP:=\{(I,J)\in\cali~|~m_\le^*\in I\}$ and $\cali^\MF:=\cali\setminus\cali^\SP$
it is clear that
\begin{eqnarray}
  x\in M(y) \mbox{ is a stationary point if and only if } \cali^\SP \mbox{ is nonempty,}\label{eq:charSP}\hfill\\
  \mbox{MFCQ is violated at } x\in M(y) \mbox{ if and only if }\cali^\MF\mbox{ is nonempty.}\label{eq:charMF}\hfill
\end{eqnarray}
(The latter characterizations may also be taken as definitions for stationary points and MFCQ.)

Second, we need a symbol for the problem size; here we simply use $(n,m)$ or $(n,m,p)$ (depending on the context),
$n$ standing for the dimension of the state variable $x$, $m=(m_\le,m_=)$ for the
number of (inequality and equality) constraints, and $p$ for the dimension of the parameter vector $y$.
For ease of notation put $m^*:=(m_\le^*,m_=)$.

Third, we need a topology on the set of problem data $(f,g,h)$ of a given problem size $(n,m,p)$
in order to specify the (large) class of problems our results refer to.
To this end we use the Whitney $C^k$ topology on $C^\infty(\R^a,\R^b)$,
where a base neighborhood of $F\in C^\infty(\R^a,\R^b)$ is defined by means of a strictly positive continuous
function $\pphi:\R^a\to (0,\infty)$ by setting
\[
U_\pphi(F):=\{G~|~|\partial_\alpha(G-F)(x)|<\pphi(x)~\forall x\in\R^a,\,|\alpha|\le k\},
\]
where $\alpha$ denotes the multi-index of a partial derivative and $|\alpha|$ its order.
Note the Whitney $C^k$ topology gets stronger (bigger), as $k$ increases.
The Whitney $C^\infty$ topology is then defined as the union of all $C^k$ topologies, $k\in\N$.

Fourth, we describe the small class of quadratic problems that turns out to already represent all
the complexity a stationary point set may have (for the problem size our theorem is referring to).
We call a problem $\calp(f,g,h)$ of type (\ref{eq:calp}) a {\em special quadratic problem} $\SQP$
if \\[.4em]
$\ds\hfill
f(x,y)=\|x-c\|_2^2,\hfill
g_i(x,y)=a_i(y)^\top x+\alpha_i(y),\hfill
h_j(x,y)=b_j(y)^\top x+\beta_j(y),\hfill
$\\[.4em]
for given fixed vector $c\in\R^n$, variable vectors $a_i(y), b_j(y)\in\R^n$
and variable scalars $\alpha_i(y),\beta_j(y)$.
We say that $\SQP$ is smooth if the data $a_i,b_j,\alpha_i,\beta_j$ are.
Sometimes it is convenient to identify $\SQP$ with its problem data $(f,g,h)$,
in particular $\SP(\SQP):=\SP(f,g,h)$ is just the stationary point set
of the special quadratic problem
$\SQP$.

Finally, we need a concept to compare the local topological structure of sets.
In order to do so we use stratifications.
A {\rm stratification} of a given set $A\subset\R^n$ is a locally finite partition $\cala\subset 2^ A$
of $A$ into smooth sub-manifolds of $\R^n$.
The latter sub-manifolds are also called {\em strata}.
In the latter case the pair $(A,\cala)$ is called a {\em stratified set}.

Two stratified sets $(A,\cala)$ and $(B,\calb)$ are called {\em isomorphic},
abbreviated by $(A,\cala)\simsim(B,\calb)$,
if there exists a homeomorphism (in general not differentiable)
$\pphi:A\to B$
such that for any $X\in\cala$ there is a $Y\in\calb$ such that
$\pphi:X\to Y$ is a diffeomorphism.
Such a homeomorphism is also referred to as an {\em isomorphism}
of the compared stratified sets, also called a {\em stratified isomorphism}.
We also write a stratified isomorphism as $\varphi:(A\cala)\to(B,\calb)$,
in order to indicate the particular stratifications to be used.
The stratified sets $(A,\cala)$ and $(B,\calb)$ are called {\em locally isomorphic}
at $a\in A$ and $b\in B$, abbreviated by $(A,\cala)_a\simsim(B,\calb)_b$,
if there exist open neighborhoods $U_a$ of $a$ and $U_b$ of $b$
and a {\em local (stratified) isomorphism}
$\pphi:(A,\cala)\cap U_a\to (B,\calb)\cap U_b$ with $\pphi(a)=b$,
where $\cala\cap U_a:=\{X\cap U_a~|~X\in\cala\}$, etc.
We do not always specify the stratifications to be used
and we then say that $A_a\simsim B_b$ are locally isomorphic if there exist stratifications $\cala$ and $\calb$
such that $(A,\cala)_a\simsim (B,\calb)_b$.

\begin{theorem}[Universality of Convex Quadratic Optimization]\label{th:main}
\hfill
\newline
For any problem size $(n,m,p)$
there exists a $C^\infty$-dense and $C^2$-open class
\goodbreak\noindent
$\calc\subset C^\infty(\R^{n+p},\R^{1+m_\le+m_=})$
such that for any choice of problem data $(f,g,h)\in\calc$ and any $(\bar x,\bar y)\in\SP(f,g,h)$
there exists a special quadratic problem $\SQP$ of the same problem size and a parameter $\bar y'\in\R^p$
such that
\[
\SP(f,g,h)_{(\bar x,\bar y)}\simsim \SP(\SQP)_{(0,\bar y')}.
\]
The local stratified isomorphism $\pphi$ can be chosen such that it preserves the
combinatorial code, i.e. $(I_0,\cali)\circ\pphi=(I_0,\cali)$,
also reading as $(I_0,\cali)(\pphi(x,y))=(I_0,\cali)(x,y)$ for all
$(x,y)\in \SP(f,g,h)_{(\bar x,\bar y)}$.
\end{theorem}

We define the {\em violation set} of MFCQ by setting
\begin{equation}\label{eq:MF}
    \MF:=\MF(g,h):=\{(x,y)~|~\mbox{MFCQ is violated at }x\in M(y)\}.
\end{equation}
Now, we consider in more detail,
to what end unstable situations compare to stable situations.
Indeed, in case that Mangasarian Fromovitz constraint qualification (MFCQ) holds
at the stationary point $0$ for $\SQP(\bar y')$,
the constraint set of $\SQP(y)$ may not become empty,
provided that $y$ is sufficiently close to $\bar y'$,
hence the continuous mapping $y\mapsto (\argmin\SQP(y),y)$ locally parameterizes $\SP(\SQP)$.
This proves that $(\SP\setminus\MF)(\SQP)$ is a topological manifold of dimension $p$.
Since the isomorphism $\pphi$ stated in our theorem is a homeomorphism preserving the
combinatorial code and since MFCQ is characterized by means of this code (recall (\ref{eq:charMF})),
Theorem \ref{th:main} implies
that $(\SP\setminus\MF)(f,g,h)$ is a topological manifold of dimension $p$ as well,
provided that the problem data are from the class $\calc$.
Note that this also holds at unstable situations of $\SP\setminus\MF(f,g,h)$.
The latter manifold property is the main result of \cite{lit:SP=manifold},
where a direct, but much more complex proof was presented.

The following even stronger topological property of stationary point sets also follows as a consequence
of Theorem \ref{th:main}, but by far not as immediately.

\begin{theorem}[Manifold with boundary property]\label{th:manifold}
For any problem size $(n,m,p)$ with $m_==0$ there exists a
$C^\infty$-dense and $C^2$-open set $\calc$ of problem data such that for any choice
of problem data $(f,g)\in\calc$ the closure $\overline\SP$ of the
stationary point set constitutes a $p$-dimensional topological manifold
with boundary. The boundary is precisely $\MF(g)$.
\end{theorem}

The proof of Theorem \ref{th:manifold} is based on ``regular transformations''
between stationary points sets on one hand and violation sets of MFCQ on the other hand.
Section \ref{sec:regtrans} is devoted to this concept.
Indeed, we will see that stationary point sets and
MFCQ violation sets are in certain sense
of the same topological complexity.
In order to explain this more exactly,
we associate the symbols $\SP$ and $\MF$ with the
corresponding problem size.

Provided that the problem data are not ill posed, we will in fact show
that any $\SP(n,m,p)$ is locally isomorphic to some $\MF(n,m^*,p+1)_{\{m_\le^*\}}$
and any $\MF(n,m,p)$ with $m_==0$ is locally isomorphic to some
$\SP(n+1,m,p-1)$,
where $\MF_I:=\{(x,y)\in\MF~|~I\subset I_0(x,y)\}$ for $I\subset\{1\dots m_\le^*\}$.

This shows that, regarding the topological complexity of $\SP(n,m,p)$,
it does not matter where a specific function enters the definition of $\SP$:
as objective function or as an inequality constraint function.
This indicates certain symmetries between objective functions and inequality constraint functions,
as known in the case of generalized semi-infinite optimization, cf. \cite{lit:gsip}.

\section{Preliminaries (from differential topology)}\label{sec:pre}
In this section we consider certain tools from differential topology
that are needed to prove our results.
In order to get an idea, to what end the techniques,
treated in this section, can be used,
please read the first paragraph of Section \ref{sec:proof1} first.

We introduce Whitney regular stratifications and fibers thereof,
which are just inverse images of single points w.r.t. local projection mappings
onto the strata.
It turns out that the topological structure of a fiber does not depend on the
particular projection mapping that is used, see Corollary \ref{cor:fiber},
and that the local structure of the entire set is just a cylinder over the fiber,
i.e. a Cartesian product of the fiber with a trivial space,
see Theorem \ref{th:thom}.
Moreover, the structure of the fibers is pulled back by transversal mappings,
see Lemma \ref{lem:pullback}.

Our main tool in this respect is Thom's Isotopy Lemma
for Whitney regularly stratified sets, c.f. \cite{lit:gwdpl}.
We say that a stratified set $(A,\cala)\subset\R^n$
(or the stratification $\cala$ of $A$) is {\em Whitney regular}
if for any pair $X,Y\in\cala$ of different strata and any point $\bar x\in X\cap\overline{Y}$
the following holds:
for any pair of sequences $x_k,y_k\to\bar x$, $x_k\in X$, $y_k\in Y$,
with the convergence properties $\span(y_k-x_k)\to L$ and $T_{y_k}Y\to T$,
we have $L\subset T$. Here, $T_yY$ stands for the tangent space of the sub-manifold $Y\subset\R^n$,
regarded as a linear subspace of $\R^n$, and we say that a sequence of $m$-dimensional subspaces $V_k\subset\R^n$
converges to an $m$-dimensional subspace $V\subset\R^n$ if the Hausdorff distance
of their intersection with the unit ball in $\R^n$ tends to zero.
(This is just the convergence in the Grassmannian manifold of $m$-dimensional
subspaces of $\R^n$.)
For a stratified set $(A,\cala)$ and a manifold $P$
we define the {Cartesian product $P\times(A,\cala):=(P\times A,P\times\cala)$},
where $P\times\cala:=\{P\times X~|~X\in\cala\}$. Note that $P\times(A,\cala)$
is again Whitney regular provided that $\cala$ is.

For splitting the local topological structure of a Whitney regularly stratified set $(A,\cala)\subset\R^n$
at $x\in X$, $X\in\cala$, into a smooth part (just $X$)
and a more complex part (for instance $A\cap T_xX^\perp$) we need tubular neighborhoods.
Given a sub-manifold $X\subset\R^n$, we call the triple
$(T,\pi,\rho)$ a tube (at $X$) if $T\supset X$ is an open neighborhood,
$\pi:T\to X$ is a projection mapping (i.e. $\pi|_X=\id_X$)
and $\rho:T\to\R$ is a quadratic distance function of a Riemannian metric (c.f. \cite{lit:gwdpl}).
For a given tube $(T,\pi,\rho)$ at $X$ and $x\in X$ let $A\cap\pi^{-1}(x)$ denote
the {\em fiber} of $A$ generated by the given tube at $x$.
In case that $X$ is a stratum of a Whitney regular stratification $\cala$ of $A$,
the stratification $\cala$ induces a natural Whitney regular stratification
$\cala\cap\pi^{-1}(x):=\{Y\cap\pi^{-1}(x)~|~Y\in\cala\}$ of the fiber $A\cap\pi^{-1}(x)$
(at least in the neighborhood
of $x$, where $(\pi,\rho):Y\to X\times (0,\infty)$
remains a submersion for any $Y\in\cala$, $Y\neq X$), c.f. {\cite{lit:jjt}}.

\begin{definition}[cf. \cite{lit:gwdpl}]\label{def:loctriv}
Let $(A,\cala)\subset U\subset\R^n$ be a stratified subset of the open set $U$
and $f:U\to P$ a mapping to a manifold $P$.
We say that $(A,\cala)$ is {\em (topologically) trivial} over $P$ (w.r.t. $f$)
if there exists a stratified set $(B,\calb)$ and a stratified isomorphism
 $\pphi:P\times(B,\calb)\to(A,\cala)$
such that $f\circ\pphi=\Pi_P$, where $\Pi_P:P\times B\to P$
denotes the canonical projection $(p,b)\mapsto p$.
The mapping $\pphi$ is referred to as a {\em trivialization}.

For any $p\in P$ the mapping $\pphi_p:B\to A$, defined by setting
$\pphi_p(b):=\pphi(p,b)$, constitutes a stratified isomorphism
between $(B,\calb)$ and $(A,\cala)\cap f^{-1}(p)$.
Therefore, for any choice of $p_1,p_2\in P$ the composition $\pphi_{p_2}\circ\pphi_{p_1}^{-1}$
is a stratified isomorphism between $(A,\cala)\cap f^{-1}(p_1)$ and $(A,\cala)\cap f^{-1}(p_2)$,
compatible with $\cala$.
Here, a stratified isomorphism between
$(A,\cala)\cap f^{-1}(p_1)$ and $(A,\cala)\cap f^{-1}(p_2)$
is called {\em compatible} with $\cala$, if for any stratum $Y$ it sends
$Y\cap f^{-1}(p_1)$ to $Y\cap f^{-1}(p_2)$.

For an open set $U\subset\R^n$ we define the restriction of $(A,\cala)$ to $U$
by setting $(A,\cala)\cap U:=(A\cap U,\cala\cap U)$, where
$\cala\cap U:=\{X\cap U~|~X\in\cala\}$.
We say that $(A,\cala)$ is {\em locally trivial} over $P$ (w.r.t. $f$)
at $a\in A$ if the point $a$ has a neighborhood $U\subset\R^n$
and $f(a)$ a neighborhood $V\subset P$ such that
$(A,\cala)\cap U$
is trivial over $V$. Here, $\pphi:V\times B\to (A,\cala)\cap U$ is called a
{\em local trivialization}.
\end{definition}

\begin{theorem}[Thom's Isotopy Lemma, cf. \cite{lit:gwdpl}]\label{th:thom}
Let $(A,\cala)$ be a locally closed Whitney regularly stratified subset
of an open set $U\subset\R^n$,
and $f:U\to P$ such that for any $X\subset\cala$ the restriction $f|_X$ is a submersion
and $f|_{\overline X\cap A}$ is a proper map, i.e. the inverse image of any compact set
in $P$ is a compact subset of $\overline X\cap A$.
Then, $(A,\cala)$ is locally trivial over $P$ (w.r.t. $f$).
\end{theorem}

The following theorem connects locally defined tubes at the same stratum
to one single tube, thus enabling Thom's Isotopy Lemma to compare their fibers
(see Corollary \ref{cor:fiber} below.)

\begin{theorem}[cf.\cite{lit:gwdpl}]\label{th:tube}
Let $X\subset U$ be a sub-manifold of an open subset $U\subset\R^n$
and $f:U\to P$ a mapping to a manifold.
Suppose that $f|_X$ is a submersion.
Now, let $X_0$, $X_1$ be relatively open subsets of $X$ such that $\overline{X}_1\subset X_0$
and $\pi_0$ is the projection mapping of a tube at $X_0$,
compatible with $f$, i.e. $f\circ\pi_0=f$.
Then, there exists a tube at $X$, also compatible with $f$, such that for its projection
mapping $\pi$ we have $\pi|_{\pi^{-1}(X_1)}=\pi_0$.
\end{theorem}

\begin{corollary}[cf. \cite{lit:gmcph}]\label{cor:tube}
Let $(A,\cala)$ be a locally closed Whitney regularly stratified subset of $\R^n$
and $(T,\pi,\rho)$ a tube at $X\in\cala$.
Then, $(A,\cala)$ is locally trivial over $X$ w.r.t. $\pi$.
(Note that $\pi$ does not need to be proper.)
\end{corollary}
\proof
(This easy proof has been taken (for the sake of more self-containment) from \cite{lit:gmcph}.)
Since, for a stratum $Y\neq X$, $(\pi,\rho)|_Y$ is a submersion in a neighborhood of $X$,
there exists for given $x\in X$ a relatively open neighborhood $X_0\subset X$ of $x$
and $\veps>0$ such that
\[
\calb:=\{X_0\}\cup\{Y\cap\rho^{-1}(0,\veps),\,Y\cap\rho^{-1}(\veps)~|~Y\in\cala,\,Y\neq X\}\cap\pi^{-1}(X_0)
\]
is a Whitney regular stratification of $B:=A\cap\pi^{-1}(X_0)\cap\rho^{-1}[0,\veps]$ and the restriction of $\pi$ to any stratum of $\calb$
is a submersion.
Restricted to $B$, $\pi$ is even a proper map, thus Theorem \ref{th:thom} applies.
\QED

\begin{corollary}\label{cor:fiber}
Let $(A,\cala)$ be a locally closed Whitney regularly stratified subset of an open set
$U\subset\R^n$.
Let $X\in\cala$ be a connected stratum and let $(T_i,\pi_i,\rho_i)$, $i=1,2$,
be tubes at relatively open subsets $X_i\subset X$.
Let $x_i\in X_i$ and
$(B_i,\calb_i):=(A,\cala)\cap\pi_i^{-1}(x_i)$ be the corresponding fibers.

Then we have $(B_1,\calb_1)_{x_1}\simsim (B_2,\calb_2)_{x_2}$,
where the local stratified isomorphism $\pphi$ can be chosen
{\em compatible} with $\cala$.
This means in particular, that for a connected stratum $X\in\cala$ the fiber
$(A,\cala)_{X,x}:=(A,\cala)\cap\pi^{-1}(x)$
generated by a tube $(T,\pi,\rho)$ at $X$ with $x\in X$ is defined up to local stratified isomorphism
without specifying the tube nor the base point $x\in X$ that is used.
\end{corollary}

\begin{remark}\label{rem:connected}
For the definition of the equivalence class $(A,\cala)_{X,x}$
it is not required that $X$ is connected, since $U\ni x$ can always be chosen
such that $X\cap U$ is connected.
\end{remark}

\noindent
{\em Proof of Corollary \ref{cor:fiber}.}
Since, for any $Y\in\cala$, $\pi_i|_Y$ is a submersion in a neighborhood of $X_i$,
$(B_i,\calb_i)$ is a Whitney regularly stratified set
in a neighborhood of $x_i$, cf. \cite{lit:jjt}.
If $X=\{x\}$ is a singleton, there exists only one single projection mapping $\pi:U\to X$,
i.e. there is nothing to prove.
Otherwise, we treat two cases.

\noindent
\underline{Case 1: $x_1\neq x_2$}.
Without loss of generality we may assume that $T_1\cap T_2=\emptyset$.
Applying Theorem \ref{th:tube} to $X_0:=X_1\cup X_2$, $\pi_0:=\pi_i$ on $\pi_i^{-1}(X_i)$
and a smaller relatively open subset of $X$ containing both $x_1$ and $x_2$,
we see that $\pi_1$ and $\pi_2$ are restrictions of a globally defined projection mapping $\pi$ of a tube at $X$
to neighborhoods of $x_1$, $x_2$, respectively,

Since $X$ is connected, there is a continuous curve $\gamma:[1,2]\to X$
with $\gamma(1)=x_1$ and $\gamma(2)=x_2$.
According to Corollary \ref{cor:tube}, any $\gamma(t)$, $t\in[1,2]$, has a relatively open neighborhood
$U_t$ in $X$ such that on $U_t$ all the fibers $(A,\cala)\cap\pi^{-1}(x)$ with a base point $x\in U_t$
are locally isomorphic and the isomorphism can always be chosen compatible with $\cala$.

Since we only need a finite number of such $U_t$ to cover $\gamma[1,2]$,
we have $(A,\cala)\cap\pi^{-1}(x_1)\simsim (A,\cala)\cap\pi^{-1}(x_2)$
and the local isomorphism may also be chosen compatible with $\cala$.

\underline{Case 2: $x_1= x_2$}.
Applying Corollary \ref{cor:tube}, we ensure the existence of $x_3\in X$, $x_3\neq x_1=x_2$,
such that
$(A,\cala)\cap\pi_1^{-1}(x_1)\simsim(A,\cala)\cap\pi_1^{-1}(x_3)$.
Now we can apply the assertion proved in Case 1 to $x_2\neq x_3$ delivering
$(A,\cala)\cap\pi_1^{-1}(x_3)\simsim(A,\cala)\cap\pi_2^{-1}(x_2)$.
Since all local isomorphisms can be chosen compatible with $\cala$, the desired isomorphism
is established.
\QED

\begin{definition}[cf. \cite{lit:jjt}]\label{def:pullback}
Let be given a stratified set $(A,\cala)\subset\R^m$ and
$f:\R^n\to\R^m$. We say that $f$ meets
a stratum $X\in\cala$ {\em transversally} (symbol $f\transversal X$) if
for any $y\in f^{-1}(X)$ we have $Df(y)[\R^n]+T_{f(y)}X=\R^m$
and we say that $f$ meets $(A,\cala)$ {\em transversally} (symbol
$f\transversal(A,\cala)$) if $f\transversal X$ for all $X\in\cala$.
\end{definition}

\begin{lemma}\label{lem:pullback}
Let be given a locally closed Whitney regularly stratified set $(A,\cala)\subset\R^m$ and
$f:\R^n\to\R^m$ such that $f\transversal(A,\cala)$.
For $X\subset\R^m$ let be $X':=f^{-1}(X)$, and define $\cala':=\{X'~|~X\in\cala\}$.
Then, $f^{-1}(A,\cala):=(A',\cala')$ is a Whitney regularly stratified set.
Moreover, for any $X\in\cala$, $x\in X$ and $x'\in X'$ with $f(x')=x$, the fibers
$(A',\cala')_{X',x'}$ and $(A,\cala)_{X,x}$ are
isomorphic and the corresponding isomorphism can be chosen {\em compatible with} $f$,
i.e. it sends any stratum $Y'\in\cala'$ to the corresponding stratum $Y\in\cala$.
\end{lemma}
\proof
The mapping $\pphi:=(id,f):\R^n\to\R^n\times\R^m$ is an embedding,
i.e. $\pphi:\R^n\to M:=\pphi(\R^n)$ is a diffeomorphism.
Defining $\chii f:\R^n\times\R^m\to\R^m$ by setting
$\chii f(y,z):=z-f(y)$,
we obviously have $M=\chii f^{-1}(0)$.
The transversality of $f$ to $Y\in\cala$
is equivalent with the property that the restriction of $\chii f$ to the stratum $\R^n\times Y$
of $\R^n\times(A,\cala)$ is a submersion, cf. \cite{lit:jjt}. This implies that
$\pphi(A',\cala')=[\R^n\times(A,\cala)]\cap\chii f^{-1}(0)$ is a Whitney regularly stratified set,
hence $(A',\cala')$ is, see \cite{lit:jjt}.
Note that $A'$ is locally closed.

The desired isomorphism between the compared fibers is established by means of a sequence
of isomorphisms. We start with a projection mapping
$\pi_1:\R^n\times\R^m\to\R^n\times X$ that is compatible with $\chii f$.
According to Theorem \ref{th:tube} such a projection mapping exists.
Since $\pi_1$ is compatible with $\chii f$ and it holds $M=\chii f^{-1}(0)$,
the restriction of $\pi_1$ to $M$ is a projection mapping onto $X\cap M$ and the
fiber
$
[\R^n\times (A,\cala)]\cap \pi_1^{-1}(x',x)
$
completely lies inside $M$.

Defining $\pi_2:\R^n\to X'$ by setting $\pi_2:=\pphi^{-1}\circ\pi_1\circ\pphi$,
it is clear that the diffeomorphism $\pphi:\R^n\to M$ transfers the representative
$(A',\cala')\cap\pi_2^{-1}(x')$ of $(A',\cala')_{X',x'}$ onto
$[\R^n\times(A,\cala)]\cap\pi_1^{-1}(x',x)$.
Note that any $Y'\in\cala'$ is transferred to the corresponding $\R^n\times Y$.

Given an arbitrary projection mapping $\pi_3:\R^m\to X$,
we know that $(A,\cala)\cap\pi_3^{-1}(x)$ represents $(A,\cala)_{X,x}$.
We define a projection mapping $\pi_4:\R^n\times\R^m\to\R^n\times X$
by setting $\pi_4(y,z):=(y,\pi_3(z))$.
Then, we obviously have
\[
[\R^n\times(A,\cala)]\cap\pi_4^{-1}(x',x) = \{x'\}\times [(A,\cala)\cap\pi_3^{-1}(x)],
\]
i.e. the smooth embedding $z\mapsto(x',z)$ is an isomorphism between
$(A,\cala)\cap\pi_3^{-1}(x)$, representing $(A,\cala)_{X,x}$,
and $[\R^n\times(A,\cala)]\cap\pi_4^{-1}(x',x)$.
Note that $Y\in\cala$ is mapped to $\R^n\times Y$.

Finally,
$[\R^n\times(A,\cala)]\cap\pi_4^{-1}(x',x)$ and $[\R^n\times(A,\cala)]\cap\pi_1^{-1}(x',x)$
are isomorphic according to Corollary \ref{cor:fiber} and Remark \ref{rem:connected},
where the isomorphism maps $\R^n\times Y$ to $\R^n\times Y$ for any $Y\in\cala$.

Altogether,
we have a sequence of isomorphisms
\begin{eqnarray*}
 (A',\cala')_{X',x'}&\to&[\R^n\times(A,\cala)]\cap\pi_1^{-1}(x',x), \\
 \mbox{}[\R^n\times(A,\cala)]\cap\pi_1^{-1}(x',x)&\to&[\R^n\times(A,\cala)]\cap\pi_4^{-1}(x',x),\\
 \mbox{}[\R^n\times(A,\cala)]\cap\pi_4^{-1}(x',x) &\to& (A,\cala)_{X,x},
\end{eqnarray*}
transferring $Y'\to\R^n\times Y\to\R^n\times Y\to Y$ for any $Y\in\cala$.
\QED

\begin{corollary}\label{cor:product}
 For any locally closed Whitney regularly stratified set $(A,\cala)\subset\R^n$,
any stratum $X\in\cala$, any base point $\bar x\in X$ and any $\bar y\in\R^m$
it holds
\[
 (A,\cala)_{X,\bar x}\simsim\left((A,\cala)\times\R^m\right)_{X\times\R^m,(\bar x,\bar y)}.
\]
\end{corollary}
\proof
Apply Lemma \ref{lem:pullback} to $f:x\mapsto(x,\bar y)$.
\QED

\section{Proof of the Universality Theorem}\label{sec:proof1}
For giving an outline,
the proof is based on a representation of $\SP$ in terms of jet extensions.
Indeed, we express $\SP$ as the inverse image of a characteristic set $\chii\SP$
w.r.t. a mapping (the jet extension) $\jsp$.
The jet extension depends on the problem data $(f,g,h)$,
whereas the characteristic set only depends on the problem size $(n,m)$.
For problem data in general position we can apply Lemma \ref{lem:pullback}
in order to verify that the local topological structure
of $\SP$ at $(\bar x,\bar y)$ is completely determined by the local
structure of $\chii\SP$ at $\bar z:=\jsp(\bar x,\bar y)$.
Having this in mind, the proof is completed by finding
an $\SQP$ whose local stationary point set is also determined
by $\chii\SP$ at $\bar z$.

In order to define $\jsp$
we collect all function values and partial derivatives of
the problem data $(f,g,h)$ that are needed to define the stationary point set
in one single mapping, called a (reduced) jet extension $\jsp:=\jsp(f,g,h)$
defined by setting
\begin{equation}\label{eq:jsp}
    \jsp:=\left((D_xg_i,g_i,D_xh_j,h_j)_{(i,j)={(1,1)\dots m}},D_xf\right),
\end{equation}
meaning that $i$ and $j$ run through all numbers from $1$ to $m_\le$ and $m_=$,
respectively.
The jet extension is called ``reduced'' since not all
partial derivatives of all functions enter.
Hence $\jsp$ maps $\R^{n+p}$ to the jet space $\JSP:=\JSP(n,m):=\R^{nm_\le+m_\le+nm_=+m_=+n}$.
Variables in the jet space are denoted as follows in order to fit easily with the mapping $\jsp$:
for $a_i,b_j\in\R^n$, $\alpha\in\R^{m_\le}$, $\beta\in\R^{m_=}$
let
\begin{equation}\label{eq:jsigma}
    \sigmasp:=
    \left((a_i,\alpha_i,b_j,\beta_j)_{(i,j)={(1,1)\dots m}},a_{m_\le^*}\right)\in\JSP.
\end{equation}
We define the combinatorial code $(I_0,\cali)$ on the elements of $\JSP$ as follows.
Let $I_0=I_0(\sigmasp):=\{i\in\{1...m\}~|~\alpha_i=0\}$, $I_0^*:=I_0\cup\{m_\le^*\}$
and $\mathop{C}(\sigmasp):=\mathop{C}$, where
\begin{equation}\label{eq:C}
    \mathop{C}:=\{(\mu,\lambda)\in[0,\infty)^{m_\le^*}\times\R^{m_=}~\left|~
    \begin{array}{l}
      \alpha\le 0,\,\mu_i=0~\forall i\not\in I_0^* \\[1mm]
      \sum\mu_i a_i+\sum\lambda_j b_j=0
    \end{array}
    \right\}.
\end{equation}
The second part $\cali=\cali(\sigmasp)$ of the code is defined to be the set of pairs $(I,J)$ with
$I\subset I_0^*$, $J\subset\{1\dots m_=\}$,
whose disjoint union $I\dot\cup J$ is minimal (w.r.t. inclusion) within the set of all
pairs such that $I\dot\cup J\neq\emptyset$ and there exists $(\mu,\lambda)\in \mathop{C}$
with $\mu_i>0~\forall i\in I$, $\lambda_j\neq 0~\forall j\in J$.
We partition $\cali$ into $\cali^\SP$ and $\cali^\MF$,
where $\cali^\SP$ contains all elements $(I,J)$ of $\cali$ with $m_\le^*\in I$.
By construction the combinatorial code is preserved by $\jsp(f,g,h)$ for all problem data
$(f,g,h)$ of size $(n,m)$, i.e. we have $(I_0,\cali)\circ\jsp(f,g,h)=(I_0,\cali)$.
In such a way we define subsets of $\JSP$ being {\em characteristic} for
stationary points and violation points of MFCQ, respectively:
\begin{eqnarray}
  \chii\SP &:=& \{\sigmasp\in\JSP~|~\cali^\SP(\sigmasp)\neq\emptyset\},\label{eq:chiSP} \\
  \chii\MF &:=& \{\sigmasp\in\JSP~|~\cali^\MF(\sigmasp)\neq\emptyset\}.\label{eq:chiMF}
\end{eqnarray}
Since $\jsp$ respects the combinatorial code,
we have $\SP=\jsp^{-1}(\chii\SP)$ and $\MF=\jsp^{-1}(\chii\MF)$.
Since the objective function $f$ does not enter the definition of $\MF$,
we also (preferably in Section \ref{sec:regtrans}) use $\jmf$,
defined as $\jsp$ with deleted entry $D_xf$,
and the corresponding $\chii\MF$ in the target space $\JMF$ of $\jmf$
as $\chii\SP$ with deleted entry $a_{m_\le^*}$.

It holds $\sigmasp\in(\chii\SP\cup\chii\MF)$ if and only if there exists $(\mu,\lambda)\in \mathop{C}(\sigmasp)$
with $\sum\mu_i+\sum|\lambda_j|=1$. Thus a compactness argument can be used
to show that $\chii\SP\cup\chii\MF$ is closed.
Analogously it follows that $\chii\MF$ is closed.
Moreover we show that $\chii\MF\subset\overline{\chii\SP}$,
hence $\chii\SP\cup\chii\MF=\overline{\chii\SP}$.
Indeed, any $\sigmasp\in\chii\MF$ can be approximated by points from $\chii\SP$.
If $\sigmasp$ does not belong to $\chii\SP$ itself,
we take some $(I,J)\in\cali(\sigmasp)$ and proceed as follows.
First we note that $m_\le^*$ does not belong to $I$.
By the minimality of $(I,J)$ it follows that the convex cone generated by
$a_i$, $i\in I$, and $b_j$, $j\in J$, is a linear subspace of $\R^n$,
which does not contain $a_{m_\le^*}$. By subtracting a small positive multiple of
$a_{m_\le^*}$ from all the vectors $a_i$, $i\in I$, $b_j$, $j\in J$,
we obtain a point from $\chii\SP$ arbitrarily close to $\sigmasp$.

We partition $\overline{\chii\SP}$ according to the combinatorial code,
i.e. two points belong to the same partition element if and only if
their combinatorial codes coincide.
The latter partition is semi-algebraic, i.e. its elements are semi-algebraic sets.
Here, a {\em semi-algebraic} subset of $\R^n$ is a finite union of elementary
semi-algebraic sets, where an {\em elementary semi-algebraic} subset of $\R^n$ is
the finite intersection of sets of the form $f^{-1}(0)$ and $g^{-1}(-\infty,0)$,
where $f,g$ are polynomials (with real coefficients).
It is immediately clear that a finite union of semi-algebraic sets is semi-algebraic.
Therefore the sets
\[
\left.\left\{(\sigmasp,\mu,\lambda)~\right|~(\mu,\lambda)\in C(\sigmasp),\,
(I_0,\cali)(\sigmasp)=\mathop{const}\right\}
\]
are semi-algebraic.
Projections of semi-algebraic sets are again semi-algebraic
according to the Tarski-Seidenberg principle,
c.f. \cite{lit:bcr},
and hence the sets of $\sigmasp$ with constant combinatorial code
are semi-algebraic as well.

A finite semi-algebraic partition (of a semi-algebraic set) admits a finite
semi-algebraic refinement that is a Whitney regular stratification, c.f. \cite{lit:bcr}.
Applying this to our partition according to the combinatorial code,
we see that there exists a finite Whitney regular stratification $\bar\cala$ of
$\overline{\chii\SP}$ such that the combinatorial code is constant
on each stratum.
Since $\chii\SP\subset\overline{\chii\SP}$ is defined by properties
of the combinatorial code and the combinatorial code is constant on any stratum of $\bar\cala$,
there exists $\cala\subset\bar\cala$ such that $\bigcup\cala=\chii\SP$.

Let $\calc:=\{(f,g,h)~|~\jsp(f,g,h)\transversal (\overline{\chii\SP},\bar\cala)\}$.
According to the Jet Transversality Theorem, c.f. \cite{lit:jjt},
the set $\calc$ is $C^\infty$-dense and $C^2$-open.
For any choice of problem data $(f,g,h)\in\calc$ we can argue as follows.
Let $\SP:=\SP(f,g,h)$, etc., and let $(\bar x,\bar y)\in\SP$ and define $\bar z:=\jsp(\bar x,\bar y)$.
Now let $Z\in\cala$ denote the stratum containing $\bar z$.
Let $(\bar F_1,\bar\calf_1):=(\overline{\chii\SP},\bar\cala)_{Z,\bar z}$ be the fiber
of $\overline{\chii\SP}$ at $Z$ with base point $\bar z$ and
$(F_1,\calf_1)=(\chii\SP,\cala)_{Z,\bar z}$ be the part of that fiber belonging to $\chii\SP$.
Let $(\bar F_2,\bar\calf_2):=(\overline\SP,\jsp^{-1}(\bar\cala))_{Z',(\bar x,\bar y)}$
be the fiber
of $\overline\SP$ at the stratum $Z':=\jsp^{-1}(Z)$
with base point $(\bar x,\bar y)$
and $(F_2,\calf_2)=(\SP,\jsp^{-1}(\cala))_{Z',(\bar x,\bar y)}$
be the part of that fiber belonging to $\SP$.
Then Lemma \ref{lem:pullback} implies that
$(\bar F_2,\bar\calf_2)_{(\bar x,\bar y)}\simsim (\bar F_1,\bar\calf_1)_{\bar z}$,
where the isomorphism $\pphi$ can be chosen compatible
with $\bar\cala$, which implies
\begin{equation}\label{eq:F2}
(F_2,\calf_2)_{(\bar x,\bar y)} \simsim (F_1,\calf_1)_{\bar z},
\end{equation}
realized by the restriction of $\pphi$.

An application of Corollary \ref{cor:tube}
yields that $\left(\SP,\jsp^{-1}(\cala)\right)$
is (at $(\bar x,\bar y))$ locally isomorphic with $Z'\times(F_2,\calf_2)$,
where the co-dimension of $Z'$ in $\R^{n+p}$ coincides with the
co-dimension of $Z$ in $\JSP$.

For fixed $\bar z\in Z$ the sets $\left(\SP,\jsp^{-1}(\cala)\right)_{(\bar x,\bar y)}$
are hence pairwise isomorphic for any choice of $(f,g,h)\in\calc$
and of $(\bar x,\bar y)$ with $\jsp(f,g,h)(\bar x,\bar y)=\bar z$.
The isomorphism may always be chosen compatible with $\cala$.
Thus our proof is completed, if we find an $\SQP\in\calc$
with $\jsp(\SQP)(0,\bar y')=\bar z$ for certain $\bar y'\in\R^p$.

To this end we first construct a special quadratic problem $\SQP_1$.
The problem size is $(n,m,p_1)$, where $p_1:=\dim\JSP-n$.
We sort the coefficients of $y\in\R^{p_1}$ as follows:
\[
y=(a_i,\alpha_i,b_j,\beta_j)_{(i,j)=(1,1)\dots m}
\]
and we define $\SQP_1:=\calp(f,g,h)$ with
\[
 f(x,y):=\frac{1}{2}\|x+\bar a_{m_\le^*}\|^2_2,
 ~\hfill g_i(x,y):=a_i^\top x+\alpha_i,
 ~\hfill h_j(x,y):=b_j^\top x + \beta_j,
\]
where $\bar a_{m_\le^*}$ is the last entry in $\bar z$.
For $\bar y':=(\bar a_i,\bar\alpha_i,\bar b_j,\bar\beta_j)_{(i,j)=(1,1)\dots m}$,
where $\bar a_i, \bar\alpha_i$ etc., are the corresponding entries of $\bar z$,
we have $\jsp_1(0,\bar y')=\bar z$, where $\jsp_1:=\jsp(\SQP_1)$.
We claim that the Jacobian $D\jsp_1(0,\bar y')$ is regular
and therefore $\SQP_1$, restricted to a small neighborhood of $(0,\bar y')$, 
belongs to $\calc$.
In order to prove this, we sort the order of partial derivatives
in the Jacobian of $\jsp_1$ as
$\left((\partial a_i,\partial\alpha_i,\partial b_j,\partial\beta_j)_{(i,j)=(1,1)\dots m},\partial x\right)$
and obtain the following Jacobian at $(0,\bar y')$
(where for each family of indices $i$ and $j$
only one representative is shown):
\begin{equation}\label{matrix}
 \begin{array}{r}
  D(D^\top_xg_i)=D(a_i)\\
  D(g_i)\\
  D(D^\top_xh_j)=D(b_j)\\
  D(h_j)\\
  D(D^\top_xf)=D(x+\bar a_{m_\le})
 \end{array}
 \left(\begin{array}{ccccc}
        Id\\
        0^\top&1&&&\bar a^\top_i\\
        &&Id\\
        &&0^\top&1&\bar b^\top_j\\
        &&&&Id
       \end{array}
 \right)
\end{equation}
Equation (\ref{matrix}) should be read as follows.
In front of each line of the matrix one can find its content,
always a full differential.
For instance, the first line is the full differential of $D^\top_xg_i=a_i$,
computed at $(0,\bar y')$.
It consists of $\partial_{a_i}a_i=Id$, $\partial_{\alpha_i}a_i=0$, etc.
The second line is the full differential of $g_i=a_i^\top x+\alpha_i$,
consisting of $\partial_{a_i}g_i=\bar x^\top=0^\top$, $\partial_{\alpha_i}g_i=1$,\dots,
$\partial_xg_i=\bar a^\top_i$.
Note that $\partial_{a_i}(g_i)$ and  $\partial_{b_j}(h_j)$ vanish because $\bar x=0$;
this fact has been emphasized by the entries zero under the main diagonal in (\ref{matrix}).
The entries of the matrix that vanish for all variables $(x,y)$ have left blank.

It is clear that (\ref{matrix}) shows a regular upper triangular matrix,
hence $\jsp_1:\R^n\times\R^{p_1}\to\JSP$ is a local diffeomorphism
in a neighborhood of $(0,\bar y')$.
This completes the proof in case $p=p_1$.

In case $p>p_1$ we use the problem $\SQP_2$, defined by setting $\SQP_2(x,y):=\SQP_1(x,y_1\dots y_{p_1})$.
It is clear that $\jsp(\SQP_2)$ meets any sub-manifold of $\JSP$ (containing $\bar z$) transversally,
which completes the proof in case $p>p_1$, too.

It remains to cover the case $p<p_1$.
To this end we consider $Z_1:=\jsp_1^{-1}(Z)$. Recall that $(0,\bar y')=\jsp_1^{-1}(\bar z)$.
It suffices to find a $p$-dimensional affine subspace $V$ of $\R^{p_1}$ containing
$\bar y'$ such that the product space $\R^n\times V$ meets $Z_1$ transversally at $(0,\bar y')$,
meaning (formally) that the embedding of that product space into $\R^n\times\R^{p_1}$ meets $Z_1$ transversally.
Indeed, since $\jsp_1$ is a diffeomorphism, the latter transversality is equivalent with
$\jsp_1(\R^n\times V)\transversal Z$ at $\bar z$.
For local coordinates $y$ in $V$ and setting $\SQP_3(y):=\SQP_1(y)$, $y\in V$,
one therefore has $\jsp(SQP_3)\transversal X$ at $\bar z$.

For finding $V$ we consider the canonical projection $\Pi:\R^n\times\R^{p_1}\to\R^{p_1}$.
It suffices to show that $\Pi(Z_1)$ is a smooth manifold and
$\Pi:Z_1\to\Pi(Z_1)$ is a diffeomorphism,
since then the existence of problem data $(f,g,h)\in\calc$
with $\jsp(f,g,h)(\R^{n+p})\cap Z\neq\emptyset$
(which holds by assumption of our theorem)
implies that the co-dimension of $Z$ in $\JSP$
does not exceed $n+p$, thus the co-dimension of $\Pi(Z_1)$ in $\R^{p_1}$
is at most $p$, implying the existence of a $p$-dimensional affine subspace $V$
through $\bar y'$ meeting $\Pi(Z_1)$ transversally at that point.

It remains to show that $\Pi:Z_1\to\Pi(Z_1)$
is a (smooth) diffeomorphism.
For later use we put this assertion in the following
slightly more general lemma.

\begin{lemma}
\label{lem:proj}
We consider a problem $\calp(f,g,h)$ of size $(n,m,p)$.
Let be $X\subset\R^n\times\R^p$, $(\bar x,\bar y)\in X$,
$U$ an open neighborhood of $(\bar x,\bar y)$,
$I\subset\{1\dots m_\le\}$ and $J\subset\{1\dots m_=\}$
such that the following conditions hold.
\begin{itemize}
\item[(i)] For any $(x,y)\in X\cap U$
the point $x$ is the unique stationary point of $\calp(y)$
with $(x,y)\in U$.
\item[(ii)] For any $(x,y)\in X\cap U$
the pair $(I^*,J)$ belongs to $\cali(x,y)$.
(Recall that $I^*=I\cup\{m_\le^*\}$
and that the latter assumption implies $I\subset I_0(x,y)$.)
\item[(iii)] The Hessian
$D^2_{xx}L_{I,J}(\bar x,\bar y,\mu,\lambda)$ is regular on the subspace of $\R^n$
consisting of all vectors which are perpendicular
to all $D_xg_i(\bar x,\bar y)$ and $D_xh_j(\bar x,\bar y)$,
$i\in I$, $j\in J$, where $\mu_i$, $\lambda_j$ are the unique multipliers
with $D_xL_{I,J}(\bar x,\bar y,\mu,\lambda)=0$.
\end{itemize}
Letting $\Pi:\R^n\times\R^p\to\R^p$ denote the canonical projection onto the second factor
and $Y:=\Pi(X)$,
there exist an open neighborhood $U'$ of $(\bar x,\bar y)$
and a smooth mapping $x:V:=\Pi(U')\to U'$
such that it holds $(x(y),y)\in X$ for any $y\in Y\cap V$.
\end{lemma}
\proof
We consider the problem
\[
\calp_{I,J}(y)\quad \mbox{minimize}~~ f(x,y) \mbox{ s.t. }g_i\le 0,\, i\in I,~h_j=0,\, j\in J.
\]
From the assumption $(I^*,J)\in\cali(x,y)$ it follows that $x$ is a stationary point also for
$\calp_{I,J}(y)$, whenever $(x,y)\in X$.
Our assumptions guarantee that the point $\bar x$
is even a {\em non-degenerate critical point} for $\calp_{I,J}(\bar y)$
in the sense of \cite{lit:jmt}.
This, according to \cite{lit:jmt}, implies
the existence of an open neighborhood $U'\subset U$ of $(\bar x,\bar y)$
and a smooth mapping $x:V:=\Pi(U')\to U'$ with the property,
that for any $y\in V$ the point $x=x(y)$ is the unique stationary point
of $\calp_{I,J}(y)$ with $(x,y)\in U'$.
Since for $(x,y)\in X\subset U'$ the point x is the only stationary point of $\calp(y)$
with $(x,y)\in U'$ and $x$ also is a stationary point of $\calp_{I,J}(y)$,
we necessarily have $x=x(y)$.This completes the proof of the lemma.\QED

With $X:=Z_1$ and $U\subset\R^n\times\R^p$ being an open neighborhood of $(0,\bar y')$,
on which $\pphi$ is a diffeomorphism,
the assumptions of the lemma are satisfied.
In fact, (ii) holds, since the combinatorial code is constant on $X$.
Since $\SQP_1(y)$ is always a quadratic
problem with strictly convex objective function,
condition (i) and (iii) are obviously guaranteed.
Since $Z_1$ is a smooth manifold, the assertion of the lemma
shows that the restriction $\Pi|_{Z_1}:Z_1\to\Pi(Z_1)$ is a smooth diffeomorphism.\QED

In our proof we also have derived the following.
\begin{corollary}
Any possible fiber $\chii\SP_{Z,z}$ of the characteristic set appears
as a fiber of $\SP(\SQP)$ a special quadratic problem $\SQP$.
\end{corollary}

\section{Regular Transformations}\label{sec:regtrans}

By means of the concept of a ``regular transformation''
we generalize and formalize what is known in optimization
theory as the ``introduction of slack variables''.

Slack variables are frequently used to transform arbitrary
finite systems of (linear) equalities and inequalities
into a standard form, where inequalities merely appear
in the form $x_i\ge 0$.
This standard form supports for instance the (academic)
implementation of the simplex method, c.f. e.g. \cite{lit:jmt}.
The ``introduction of slack variables''
roughly works as follows. Given a function $f:\R^n\to\R$
we are considering the set
\begin{equation}
\label{eq:slack1}
\{x\in\R^n~|~f(x)\ge 0\}.
\end{equation}
Using a slack variable $x_{n+1}\in\R$ we may define
the functions $f_1'(x,x_{n+1}):=x_{n+1}$ and
$f_2'(x,x_{n+1}):=x_{n+1}-f(x)$ and a transformed set
\begin{equation}\label{eq:slack2}
\{(x,x_{n+1})\in\R^n\times\R~|~f_2'(x,x_{n+1})=0,~f_1'(x,x_{n+1})\ge 0\}.
\end{equation}
It is easy to see that $x\in\R^n$ belongs to the original set in (\ref{eq:slack1})
if and only if the pair $(x,x_{n+1}):=(x,f(x))$ belongs to the transformed set
in (\ref{eq:slack2}) and that the transforming mapping
$x\mapsto(x,x_{n+1})=(x,f(x))$ is as smooth as $f$.
Altogether, the introduction of slack variables turns sets given
in the form (\ref{eq:slack1}) into sets in the form (\ref{eq:slack2}).

Generalizing this we will turn stationary points of a given problem $\calp(f,g,h)$
into points where MFCQ is violated by just replacing the objective function
$f(x,y)$ by a new constraint function $g_{m_\le^*}(x,y,y_{p+1}):=f(x,y)-y_{p+1}$.
Moreover we will turn violation points of MFCQ of $M(g)$ into
stationary points with MFCQ by replacing the constraints $g_i(x,y)\le 0$
by $g_i(x,y)-x_{n+1}\le 0$ and adding an objective function $f(x,x_{n+1},y):=x_{n+1}$
to be minimized.

What is common in these transformations? In any case,
a new variable is introduced, smoothly depending on the original ones.
In the slack variable case we already described this mapping as $x\mapsto x_{n+1}:=f(x)$.
In the other two examples the corresponding mappings are
$(x,y)\mapsto y_{p+1}:=f(x,y)$ and $(x,y)\mapsto x_{n+1}:=0$,
respectively.

This gives rise to the concept of a ``regular transformation'',
defined as follows.
The starting point of the transformation is the object that is transformed.
It consists of a fixed {\em characteristic set} $F\subset\R^m$
in the target space of {\em defining mappings} $f:\R^n\to\R^m$,
which still can be chosen from a function space.
The characteristic set is supposed to admit a Whitney regular stratification.

In applications, $F$ may be the characteristic set $\chii\SP(n,m)$
for stationary points of a given problem size, see (\ref{eq:chiSP}),
and $f$ the reduced jet extension $\jsp(n,m,p)$,
defining the stationary point set corresponding to a particular choice
of problem data.
Since $f^{-1}(F)$, by Lemma \ref{lem:pullback},
inherits its local topological structure from
the fibers of $F$,
our main concern are these fibers and how they transform.

Now let us formally describe,
how such transformations should work.
Part of the transformation are a smooth embedding $\Ttar:\R^m\hookrightarrow M$,
where $M\subset\R^{m'}$ is a smooth sub-manifold,
and a new characteristic set $F':=\Tset(F)\subset\R^{m'}$, such that
\begin{equation}
\label{eq:1}
\Ttar(F)=F'\cap M.
\end{equation}
This describes the action of the transformation on the target space (of $f$).
The transformation also acts on the (variable) space $\R^n$ and involves
$\Tvar:=\Tvar(f):\R^n\to\R^{n'}$
and creates a new defining function
$f':=\Tfunc(f):\R^{n'}\to\R^{m'}$
such that the following conditions hold:
\begin{equation}
\label{eq:2}
f'\circ\Tvar=\Ttar\circ f,
\end{equation}
\begin{equation}
\label{eq:3}
f'(\R^{n'})\subset M,
\end{equation}
\begin{equation}
\label{eq:4}
f'^{-1}(M')\subset N',
\end{equation}
where $M':=\Ttar(\R^m)$, $N':=\Tvar(\R^n)$.
The transformation is called {\em regular} if for any Whitney regular
stratification $\calf$ of $F$ there exists a Whitney regular stratification $\calf'$
of $F'$ such that:
\begin{equation}
\label{eq:r1}
\Ttar(F,\calf)\mbox{ coincides with }(F',\calf')\cap M,
\end{equation}
\begin{equation}
\label{eq:r2}
M\transversal(F',\calf')\mbox{ in }\R^{m'},
\end{equation}
\begin{equation}
\label{eq:r3}
f'\transversal M'\mbox{ in }M.
\end{equation}
Although seemingly quite involved,
regular transformations are in some sense natural
and not very hard to cope with.
Indeed, noting that $\Tvar:\R^n\to N'$
and $\Ttar:\R^m\to M'$
are diffeomorphisms,
we have the following commutative diagram
\begin{equation}
\nonumber
\label{eq:diagram1}
\begin{array}{ccccc}
\R^n&\overset{f}{\rightarrow}&\R^m&\hookleftarrow&F\\
{\Tvar}\updownarrow\phantom{\Tvar}&&\phantom{\Tvar}\updownarrow\Ttar&&\phantom{\Tvar}\updownarrow\Ttar\\
N'&\overset{f'}{\rightarrow}&M'&\hookleftarrow&F'\cap M'.
\end{array}
\end{equation}
Conditions (\ref{eq:1})-(\ref{eq:r1}) guarantee that
$\ds f'^{-1}(F')\overset{f'}{\rightarrow}F'\cap M'$
is just
$\ds f^{-1}(F)\overset{f}{\rightarrow}F$
in other coordinates, see below for details.
Since we want to have the freedom to add additional $n'-n$ variables
(e.g. slack variables) to the variable space
we need to add additional $n'-n$ dimensions to the target space,
forming $M$,
but without enlarging $F'$,
just in order to guarantee that (as a stratified space)
the set $f'^{-1}(F')$ has the same dimension as $f^{-1}(F)$.
This is formalized in (\ref{eq:r3}).
In case of $\MFSP$,
transforming MFCQ violation point sets in $n$ dimensions
into stationary points sets in $n+1$ dimensions,
we need to add additional dimensions to the target space,
for instance in order to model the $n+1$st partial derivatives.
In order to keep transversality and the dimension of
$f'^{-1}(F')$ simultaneously,
this requires the addition of the same number of degrees of
freedom to the characteristic set $F'$.
This is formalized in (\ref{eq:r2}).

\begin{lemma}
\label{lem:compo}
Compositions of regular transformations are regular transformations.
\end{lemma}
The proof of the lemma is easy and therefore omitted.

\begin{theorem}[On regular transformations]
\label{theo:regTrafo}
Let be given a locally closed Whitney regularly stratified set
$(F,\calf)\subset\R^m$ and a smooth mapping $f:\R^n\to\R^m$ with
$f\transversal(F,\calf)$. Let $T$ be a regular transformation as described above.
Then the following assertions hold.
\begin{itemize}
\item[i)]
The diffeomorphism $\Tvar:\R^n\to N'$ is a stratified isomorphism between
$f^{-1}(F,\calf)$ and $f'^{-1}(F',\calf')$.
Equation (\ref{eq:diagram}) shows a commutative diagram
of stratified mappings.
\begin{equation}
\label{eq:diagram}
\begin{array}{ccc}
f^{-1}(F,\calf)&\overset{f}{\rightarrow}&(F,\calf)\\
\phantom{\Tvar}\downarrow\Tvar&&\phantom{\Tvar}\downarrow\Ttar\\
f'^{-1}(F',\calf')&\overset{f'}{\rightarrow}&(F',\calf')
\end{array}
\end{equation}
\item[ii)]
For $\Ttar(z)=z'$, $Z\in\calf$ and $Z'\in\calf'$ with $z\in Z$ and $z'\in Z'$,
the mapping $\Ttar$ is a stratified isomorphism between appropriate representatives
of $(F,\calf)_{Z,z}$ and $(F',\calf')_{Z',z'}$.
\end{itemize}
\end{theorem}
\proof {\em ad(i)}
First we prove that $\Tvar$ maps $f^{-1}(F)$ surjectively (thus bijectively) onto $f'^{-1}(F')$.
Indeed, since $\Ttar:\R^m\to M'$ is a diffeomorphism
and since $M'\supset\Ttar(F)\stackrel{(\ref{eq:1})}{=}F'\cap M$
we also have $\Ttar(F)=F'\cap M'$ and therefore:
\begin{equation}
\label{eq:Ttar}
x\in f^{-1}(F)\Leftrightarrow f(x)\in F\Leftrightarrow
f'\circ\Tvar(x)\stackrel{(\ref{eq:2})}{=}\Ttar\circ f(x)\in F'\cap M'.
\end{equation}
Hence $x\in f^{-1}(F)$ implies $\Tvar(x)\in f'^{-1}(F')$.
On the other hand let $x'\in f'^{-1}(F')$.
Now we prove the existence of $x\in f^{-1}(F)$
with $\Tvar(x)=x'$,
i.e. that $\Tvar:f^{-1}(F)\to f'^{-1}(F')$ is surjective as well.
By definition we have $f'(x')\in F'$,
which by (\ref{eq:3}) implies
$\ds f(x')\in F'\cap M\stackrel{(\ref{eq:1})}{=}F'\cap M'$.
Since $f(x')\in M'$, (\ref{eq:4}) yields $x'\in N'$, thus $x'=\Tvar(x)$ for
some $x\in\R^n$.
Applying (\ref{eq:Ttar}) we conclude $x\in f^{-1}(F)$.

The smooth embedding $\Ttar$, applied to $(F,\calf)$, is a {\em stratified mapping}
due to (\ref{eq:r1}),
i.e. a mapping that sends strata into strata.
The mappings $f$ and $f'$ are stratified mapping by definition.
Moreover, for any $Z\in\calf$ the mapping $\Tvar$ maps the stratum $f^{-1}(Z)$
to the stratum $f'^{-1}(Z')$, where $Z'\in\calf'$ is the stratum
containing $\Ttar(Z)$, which exists due to (\ref{eq:r1}).
Hence, by (\ref{eq:2}), also $\Tvar$ is a stratified mapping, even a stratified isomorphism,
since it is smooth and maps $f^{-1}(F)$ bijectively onto $f'^{-1}(F')$.

{\em ad (ii)}
For ease of notation we identify one point sets with its point.
$\Ttar(Z)$ is a sub-manifold of $M'\subset M\subset\R^{m'}$.
Since we have $Z'\cap M\stackrel{(\ref{eq:1})}{=}Z'\cap M'\stackrel{(\ref{eq:r1})}{=}\Ttar(Z)$
and $M\transversal Z'$ (by (\ref{eq:r2})) we can choose local coordinates such that
$z'=0$, $\Ttar(Z)=\R^{d_1}\times 0$,
$M'=\R^{d_1}\times\R^{d_2}\times 0$,
$M=\R^{d_1}\times\R^{d_2}\times\R^{d_3}\times 0$,
$\R^{m'}=\R^{d_1}\times\R^{d_2}\times\R^{d_3}\times\R^{d_4}$
and $Z'=\R^{d_1}\times 0\times 0\times\R^{d_4}$.
Let $\Pi_i:\R^{m'}\to\R^{d_i}$ denote the canonical projection onto the $i$'th factor.
For generating the fibers to be compared,
we have the freedom to choose the local projection mappings $\pi$ and $\pi'$ on the strata
$Z$ and $Z'$, respectively.
For ease of argumentation we choose
$\ds \pi:=\Ttar^{-1}\circ\Pi_1\circ \Ttar$,
i.e. it holds
\begin{equation}
\label{eq:pi}
\Ttar\circ\pi\circ\Ttar^{-1}=\Pi_1,
\end{equation}
and we choose $\pi':\R^{m'}\to Z'$ as $\pi':=(\Pi_1,0,0,\Pi_4)$.
Then we have $\Ttar(\pi^{-1}(z))=0\times\R^{d_2}\times 0\times 0$
and $\pi'^{-1}(z')=0\times\R^{d_2}\times\R^{d_3}\times 0$.
By (\ref{eq:1}) it holds $F'\cap M\subset M'$, thus
$\pi'^{-1}(z')\cap F'=(0\times\R^{d_2}\times 0\times 0)\cap F'$.
Due to (\ref{eq:1}) and (\ref{eq:pi}) we also have
$\Ttar(\pi^{-1}(z)\cap F)=(\Pi_1|_{M'}^{-1}(0)\times 0\times 0)\cap F'
=(0\times\R^{d_2}\times 0\times 0)\cap F'$.
Altogether we see that $\Ttar$ maps $(F,\calf)\cap\pi^{-1}(x)$
onto $(F',\calf')\cap\pi'^{-1}(z')$,
proving the assertion.
\QED

In the following we treat transformations that we call ``elementary''.
We do not specify all things explicitly, since some of them
are given by the concept of a transformation itself, for instance
$M'\subset\R^{m'}$ is already given if $\Ttar$ is defined.
If not defined differently (Type 2), $M$ coincides with $\R^{m'}$.
We always consider the situation of an implicitly defined (family)
$f^{-1}(F)$ with fixed characteristic set $F$ in the target space of $f$.

\begin{definition}[Elementary Transformation of Type 1]
\label{def:ET1}
Type 1 transformations introduce a new variable. They are induced by a
smooth function $g:\R^n\to\R^p$ as follows.
We define $\Tvar(x):=(x,g(x))\in\R^{n'}$, $n':=n+p$.
In the target space this corresponds $\Ttar(z):=(z,0)\in\R^{m'}$,
$m':=m+p$.
In order to have (\ref{eq:1}) we need to define $F':=\Ttar(F)=F\times 0$.
and for (\ref{eq:r1}) let $\calf':=\Ttar(\calf):=\{\Ttar(Z)~|~Z\in\calf\}$
for given stratification $\calf$ of $F$.
For the definition of $f'$ we distinguish two subtypes.
In sub-case 1a the latter function is defined by setting
$f'(x,y):=(f(x),\pm(f(x)-y))$
and in sub-case 1b by
$f'(x,y):=(y,\pm(f(x)-y))$.
\end{definition}

\begin{definition}[Elementary Transformation of Type 2]
\label{def:ET2}
Type 2 transformations just expand the target space by some additional dimensions.
They are defined by setting $m':=m+q$ and for given $Q\in\R^q$
by $\Ttar(z):=(z,Q)$.
We put $M=M':=\R^m\times Q$, $F':=F\times \R^q$
and for given stratification $\calf$ of $F$
we put $\calf':=\calf\times\R^q$.
Finally we define $\Tvar(x):=\id_{\R^n}$ and $f'(x):=(x,Q)$.
\end{definition}

\begin{definition}[Elementary Transformation of Type 3]
\label{def:ET3}
Type 3 transformations deform the target space by means of a diffeomorphism.
In order to define them let $U\subset\R^m$ be an open set such that $U\supset F\cap f(\R^n)$
and $\pphi:U\to V\subset\R^m$ be a diffeomorphism.
Then we define $M:=M':=\R^{m'}:=\R^m$, $\Ttar:=\pphi$,
$(F',\calf'):=\pphi(F,\calf)$ and $f':=\pphi\circ f$.
Letting $\Tvar:=\id_{\R^n}$ we are done.
\end{definition}

\begin{lemma}
\label{lem:types}
All types of elementary transformations are regular.
\end{lemma}
The proof of the lemma is straightforward and therefore omitted.

Now we show that the introduction of slack variables and the transformations
between stationary point sets and violation sets of MFCQ
can be written as (compositions of) elementary transformations
and that they are therefore regular.
In order to do so treat the respective situation as
some object to be transformed, i.e. as some $f^{-1}(F)$.
Then we describe, how the transformation acts on the
variable and the target space and also which
new defining function $f'=\Tfunc(f)$ is used.
In case of transformations between stationary point sets
and violation sets of MFCQ we are particularly interested,
how the transformation changes the combinatorial code,
since this offers the opportunity
to track specific parts of the transformed sets, which are defined by means of this code.
The latter philosophy forms the basis of the proof of Theorem \ref{th:manifold}.

The transformation $\SLACK$, introducing a slack variable,
deals with $f^{-1}(F)$,
where $f:\R^n\to\R$ and $F:=[0,\infty)$.
We define $\SLACK^{\mathop{var}}(x):=(x,f(x))$,
$\SLACK^{\mathop{tar}}(z):=(z,0)\in M:=\R^{m'}$, $m':=2$.
Setting $F':=\SLACK^{\mathop{tar}}(F)=[0,\infty)\times 0$ and $f'(x,y):=(y,y-f(x))$
we see that $\SLACK$ is of type 1b.

The transformation $\SPMF$, transforming $\overline\SP(n,m,p)\subset\R^{n+p}$
into $\MF(n,m^*,p+1)$,
can be written as a type 1a transformation.
For problem data $(f,g,h)$ of $\calp(n,m,p)$ the defining mapping
for $\SP$ and $\overline\SP$ is $\jsp:=\jsp(f,g,h)$.
Then the setting to be transformed is $\jsp^{-1}(\overline{\chii\SP})$,
where $\overline{\chii\SP}:=\overline{\chii\SP}(n,m)$.
The transformation $\SPMF$ is defined to be the type 1a transformation induced
by the mapping $g_{m^*_\le}=f:\R^{n+p}\to\R$.
The target space of the mapping $\SPMFtar:\JSP=\JSP(n,m)\to\JSP\times\R$,
$\sigmasp\mapsto(\sigmasp,0)$,
may be considered as {$\JMF=\JMF(n,m^*)$}.
By the very construction it holds $(I_0,\cali)(\sigmasp,0)=(I_0^*,\cali)(\sigmasp)$,
where $I^*:=I\cup\{m^*\}$.
Thus the transformation acts uniquely on the combinatorial code,
written as $\SPMFcode((I_0,\cali)):=(I_0^*,\cali)$.
Moreover $\SPMFtar(\overline{\chii\SP})$ precisely coincides with
$\chii\MF_{\{m^*_\le\}}=\{\sigmamf\in\chii\MF~|~m^*_\le\in I_0(\sigmamf)\}$
and
{$\SPMFfunc(\jsp(f,g,h))=\jmf(G,H)$},
where $H_j,G_i:\R^n\times\R^{p}\times\R\to\R$ are the following functions:
\begin{equation}
\label{eq:SP2MF}
\begin{array}{rcll}
H_j(x,y,y_{p+1})&:=&h_j(x,y),\quad &j=1,\dots,m_=,\\
G_i(x,y,y_{p+1})&:=&g_i(x,y),\quad &i=1,\dots,m_\le,\\
G_{m_\le^*}(x,y,y_{p+1})&:=&f(x,y)-y_{p+1}.
\end{array}
\end{equation}
Therefore $\SPMFvar(\overline\SP)=\MF_{\{m^*_\le\}}=\jmf^{-1}(\chii\MF_{\{m_\le^*\}})$.

The transformation $\MFSP$ models violation sets of MFCQ
in terms of stationary point sets, provided that no equality
constraints are present. In the latter case we identify $m\equiv m_\le$.
For problem data $g$ for {$M(n,m,p)$} (with $m_==0$)
we define new problem data $(F,G)$ for $\calp(n+1,m,p)$,
acting on variables $(x,x_{n+1},y)\in\R^n\times\R\times\R^p$,
by setting
\begin{equation}
\label{eq:MFSP}
F:=x_{n+1},~G_i:=g_i(x,y)-x_{n+1},~i=1\dots m.
\end{equation}
$\MFSP$ is the composition of three elementary transformations.
Its variable space action is $\MFSPvar:\R^n\times\R^p\to\R^{n+1}\times\R^p$
defined by setting $\MFSPvar(x,y):=(x,0,y)$.
Letting $\SP:=\SP(F,G)$ and $\SP_0:=\SP\cap F^{-1}(0)$,
we obviously have $(x,y)\in\MF:=\MF(g)$
if and only if $(x,0,y)\in\SP_0$.
(In contrast to $\SPMF$ the newly introduced variable
does not contribute to the parameter space but to the state space.)
In order to express $\MFSP$ as a composition of elementary transformations
we use transformations of types 1a, 2 and 3.
For convenience (in a later application) we start with
transformation $T_2$ of type 2 followed by $T_3$ of type 3
which are known to only act on the target space
and to leave the variable $(x,y)\in\R^{n+p}$ unchanged.

We define $\TTtar:\JMF:=\JMF(n,m)\to\JSP:=\JSP(n+1,m)$
by setting
\begin{equation}
\label{eq:T2}
\TTtar:(a_i,\alpha_i)_{i=1\dots m}\mapsto
  \left(\left(
           \binom{a_i}{-1},\alpha_i
        \right)_{i=1\dots m},
        \binom{0}{1}
  \right)=:\sigmasp.
\end{equation}
Reordering variables we identify $\JSP\equiv\JMF\times\R^q$, $q:=m+n+1$,
$\sigmasp\equiv(\sigmamf,\bar Q)$,
where
\begin{eqnarray}
\label{eq:Q}
\nonumber
\bar Q&&:=\left((\bar a_{i,n+1})_{i=1\dots m},\bar a_{m^*,n+1},\bar a_{m^*}\right)\\
      &&:=\left((-1)_{i=1\dots m},1,0\right)\in\R^{m}\times\R\times\R^n=\R^q.
\end{eqnarray}
By construction it holds $\TTfunc(\jmf)(x,y)=\TTtar\circ\jmf(x,y)=\jsp(x,0,y)$,
where $\jsp:=\jsp(F,G)$.
Recall that we have $M_2'=M_2=\JMF\times \bar Q$.
By inspection we see that $T_2$ acts as follows on the combinatorial code:
$\TTcode:(I_0,\cali)\mapsto(I_0,\cali^*)$,
where $\cali^*:=\{I\cup\{m^*\}~|~I\in\cali\}$.
Since $\TTtar:\JMF\to M_2'$ is a global diffeomorphism,
this proves $\TTtar(\chii\MF)=\chii\SP\cap M_2'$.
Moreover we see that
$\TTtar(\chii\MF)\cap\chii\MF(n+1,n)=\emptyset$.
Since $T_2$ is a regular transformation we automatically have
$M\transversal\left(\TTtar(\chii\MF,\cala)\times\R^q\right)$
for any Whitney regular stratification $\cala$ of $\chii\MF$,
but we cannot automatically conclude $M\transversal(\chii\SP,\cala')$
for some Whitney regular stratification $\cala'$ of $\chii\SP$.

The latter (technical) problem is solved by means of the transformation $T_3$,
which deforms $\JSP$ locally but leaves $\chii\SP\cap M_2'$ unchanged.
Let us write a point $Q\in\R^q=\R^{m+1}\times\R^n$ as
$Q=\left((a_{i,n+1})_{i=1\dots m^*},a_{m^*}\right)$.
We define the
open subset $U:=\{Q\in\R^q~|~a_{i,n+1}<0,~i=1\dots m,~a_{m^*,n+1}>0\}\subset\R^q$
which obviously contains $\bar Q$.
Let $\pphi:\JMF\times U\to\JSP$ be defined by setting
\begin{equation}
\label{eq:pphi}
\pphi(\sigmamf,Q):=\left(
                      \left(a_{i,n+1}\binom{-a_i+v}{1},\alpha_i
                      \right)_{i=1\dots m},
                      a_{m^*,n+1}\binom{v}{1}
                   \right).
\end{equation}
One easily checks that $\pphi(.,\bar Q)=\TTtar$ and that $\pphi$ is a smooth diffeomorphism
onto its image $\pphi(\JMF\times U)=:V$.
By construction it holds
$(I_0,\cali)\left(\pphi(\sigmamf,Q)\right)=\TTcode\left((I_0,\cali)(\sigmamf)\right)$.
Hence we have $\pphi(\chii\MF\times U)=\chii\SP\cap V$.
Identifying $\JSP\equiv\JMF\times\R^q$,
we may consider $\JMF\times U$ as an open subset of $\JSP$,
thus $\pphi$ as a diffeomorphism, defined on an open set containing $M_2$,
such that $\pphi|_{M_2}=\id_{M_2}$.

Let $T_3$ be the type 3 transformation induced by $\pphi$.
Since $\TTTfunc=\id$ and $\TTTcode=\id$ it holds
$\TTTTfunc=\TTfunc$ and $\TTTTcode=\TTcode$.

In order to describe $T_1$, necessary to understand $\MFSP$
as a composition $T_1\circ T_3\circ T_2$ of elementary transformations,
we do not write $\SP$ in our standard form $\SP=\jsp^{-1}(\chii\SP)$.
Instead we use a (slightly less reduced) jet extension of $(F,G)$,
namely $\jsp_0:=(\jsp_0(F,G)):=(\jsp(F,G),F):\R^{n+1}\times\R^p\to\JSP_0:=\JSP\times\R$.
The corresponding characteristic set is $\chii\SP_0:=\chii\SP\times\R$.
Then we obviously have $\SP=\jsp_0^{-1}(\chii\SP_0)$.
Since $\MFSPvar(x,y)=(x,0,y)$ is already given,
we have $\tvar(x,y)=(x,0,y)$, thus our type 1a transformation is induced
by the mapping $\R^n\times\R^p\to\R$, $(x,y)\mapsto 0$.
Using the ``$+$'' version of $\tfunc$ we conclude that
$\MFSPfunc(\jmf(g))=\jsp_0(F,G)$.
Moreover we have $\MFSPtar(\chii\MF)=\ttar(\chii\SP\cap V)=\chii\SP_0\cap(V\times\R)$.

\begin{theorem}
\label{theo:SPMFSP}
\begin{itemize}
\item[(i)]
The transformations $\SPMF$ and $\MFSP$ (the latter only defined for $m_==0$)
are regular. The transformations are well-defined for arbitrary (smooth) problem data.
The mappings $\SPMFvar$ and $\MFSPvar$ are smooth embeddings,
mapping $\SP(n,m,p)(f,g,h)$ onto $\MF(n,m^*,p+1)(G,H)_{\{m^*_\le\}}$
and $\MF(n,m,p)(g)$ onto $\SP(n+1,m,p)(F,G)_0$, respectively.
\item[(ii)]
For any problem size $(n,m,p)$ with $m_==0$, any Whitney regular stratification
$\cala$ of $\chii\MF(n,m)$, any choice of problem data
$g$ for $M(n,m,p)$ with $\jmf(g)\transversal\cala$ and any $(\bar x,\bar y)\in\MF(g)$
the following holds.
There exist problem data $(f',g')$ for $\calp(n+1,m,p-1)$ and $y'\in\R^{p-1}$
such that $\jsp(f',g')\transversal\cala'$ for certain Whitney regular stratification
of $\chii\SP(n+1,m)$ and there exists a local stratified isomorphism
between $\MF(g)_{(\bar x,\bar y)}$ and $\SP(f',g')_{(0,y')}$.
\end{itemize}
\end{theorem}
\proof
Assertion (i) is clear from construction. In order to verify (ii),
put $\bar z:=\jmf(g)(\bar x,\bar y)$ and let $Z\in\cala$
be the stratum containing $\bar z$.
Let $d$ denote the co-dimension of $Z$ in $\JMF$.
Then the dimension of $\jmf^{-1}(Z)$ is $n+p-d$.

Let $z':=(T_3\circ T_2)^{\mathop{tar}}(\bar z)$ and $Z'\in\cala'$
be the stratum of the corresponding stratification $\cala'$
of $\chii\SP(n+1,m)$ that contains $z'$. Note that the co-dimension
of $Z'$ in $\JSP(n+1,m)$ is also $d$.
As in the proof of Theorem \ref{th:main}
we construct the problem $\SQP_1$, such that $\jsp_1\transversal\cala'$
and such the $z'=\jsp_1(0,\bar y')$ for certain $\bar y'\in\R^{\bar p}$, $\bar p:=\dim\JSP-(n+1)$.
In analogy to the proof of Theorem \ref{th:main}, for $p'$ large enough,
we can construct a $p'$-dimensional affine subspace of $\R^{\bar p}$ containing $\bar y'$
such that the restriction of $\SQP_1$ to $V$ also meets $\cala'$ transversally.
In order to obtain the same local structure as in $\SP(n,m,p)$
we only need to guarantee that the dimension of $\jsp_1|_V^{-1}(Z')=(n+1)+p'-d$ also is $n+p-d$,
being satisfied if $p'=p-1$.\QED

\section{Proof of Theorem \ref{th:manifold}}\label{sec:proof mani}

In the Introduction we have already seen that for problem data from $\calc$
the set $\SP\setminus\MF$ is a topological manifold of dimension $p$.
Since $\overline\SP\setminus\MF=\SP\setminus\MF$, it therefore
suffices to locally prove the assertion at the points from $\MF\subset\overline\SP$.

In view of Corollary \ref{cor:tube} and Lemma \ref{lem:pullback}
it suffices to prove for any Whitney regular stratification $\cala$
of $\overline{\chii\SP}$ respecting the combinatorial code,
any stratum $Z\in\cala$ with $Z\subset\MF$ and any point $z\in Z$,
that
the fiber $\overline{\chii\SP}_{Z,z}$ constitutes
a topological manifold with boundary
and that the boundary exactly is $\chii\MF_{Z,z}$.

The key of the proof consist of performing the transformation
$T=T_3\circ T_2\circ\SPMF$ and to investigate, applying Theorem \ref{theo:regTrafo} (ii),
onto which part of a fiber the target space action $\Ttar$ transforms the fiber $\chii\MF_{Z,z}$.
In this investigation our information on $\Tcode$ is crucial.
In fact we have $\Tcode=\TTTcode\circ\TTcode\circ\SPMFcode$.
$\SPMFcode$ adds $m^*$ to $I_0$ and leaves $\cali$ unchanged.
In the next step, $\TTcode$ leaves $I_0$ unchanged and adds $(m^*)^*$
to any member of $\cali$, and finally $\TTTcode$ does not change anything.

In other words,
$z'\in\Ttar(\JSP(n,m))$ belongs to the image of $\overline{\chii\SP}(n,m)$
under $\Ttar$ if and only if
the constraint with index $m^*$,
the one derived from the original objective function,
is active.
The point $z'$ even belongs to image of $\MF(n,m)$
if and only if, beside that the constraint with index $m^*$ is active,
the gradient of this constraint is not required
for making $z'$ a stationary point for $\SP(n+1,m^*)$.
This is due to the fact that (the definition of) MFCQ does not care of $f$.

After renumbering (decreasing $m$ by $1$)
it therefore suffices to prove the following lemma.

\begin{lemma}
\label{lem:prove}
For any problem size $(n,m)$ with $m_==0$,
any Whitney regular stratification $\cala$ of
$\chii\SP=\chii\SP(n,m)$ respecting the combinatorial code,
any stratum $Z\in\cala$ with $Z\cap\chii\MF=\emptyset$
and $Z\subset\chii\SP^*_{\{m\}}$,
and any point $\bar z\in Z$,
the part of the fiber $(\chii\SP,\cala)_{Z,\bar z}$
belonging to $\chii\SP_{\{m\}}$
is a topological manifold with boundary
and the latter boundary is exactly
the part of the fiber belonging to $\chii\SP^*_{\{m\}}$.
\end{lemma}
Here $z$ belongs to $\chii\SP_{\{m\}}$ if and only if
$m\in I_0(z)$ and to $\chii\SP^*_{\{m\}}$
if and only if $m\in I_0(z)$ and there exists $I\in\cali(z)$
with $m\not\in I$.

Note that the assumption $m_==0$ has merely been
made to make it easier to see, that it exactly deals with
the result of $\Ttar$.
The assertion of the lemma itself does not require this assumption.

\smallskip
\proof
The proof of the lemma is mainly based on the simple observation,
that $(z,a_m,\alpha_m)$ belongs to $\chii\SP(n,m)^*_{\{m\}}$
if and only if $z$ belongs to $\chii\SP(n,m-1)$ and $\alpha_m=0$.
This gives rise to a manifold (with boundary) chart of
$\chii\SP(n,m)$, using a special quadratic problem.

As before we write elements of $\JSP=\JSP(n,m)$
as $z=(y,a_{m^*})=\left((a_i,\alpha_i)_{i=1\dots m},a_{m^*}\right)$
and in particular
$\bar z=(\bar y,\bar a_{m^*})=\left((\bar a_i,\bar\alpha_i)_{i=1\dots m},\bar a_{m^*}\right)$.

As in the proof of Theorem \ref{th:main} we use
$\SQP=\SQP(n,m,p)$ with $p:=\dim\JSP-n$,
whose problem data $f,g_i:\R^n\times\R^p\to\R$
are defined by setting
\begin{equation}
\label{eq:sqp}
f(x,y):=\frac{1}{2}\|x-\bar a_{m^*}\|_2^2,~g_i(x,y):=x^\top a_i+\alpha_i,~i=1\dots m.
\end{equation}
Letting $\jsp:=\jsp(\SQP)$ we have $\jsp(0,\bar y)=\bar z$
and
there exist open neighborhoods $U$ of $(0,\bar y)$ and $V$ of $\bar z$
such that $\jsp:U\to V$ is a diffeomorphism.
Put $\tilde Z:=\jsp^{-1}(Z)$, $\SP:=\jsp^{-1}(\chii\SP)$,
$\SP_{\{m\}}:=\jsp^{-1}(\chii\SP_{\{m\}})$ and
$\SP^*_{\{m\}}:=\jsp^{-1}(\chii\SP^*_{\{m\}})$.
Since $\jsp$ respects the combinatorial code
it suffices to show the existence of a projection mapping $\tilde\pi:\R^n\times\R^p\to\tilde Z$,
defined on a neighborhood of $(0,\bar y)\in\tilde Z$, such that
the fiber $\SP_{\{m\}}\cap\tilde\pi^{-1}(0,\bar y)$
constitutes a topological manifold with boundary $\SP^*_{\{m\}}\cap\tilde\pi^{-1}(0,\bar y)$.
Recall that for $(x,y)\in\R^n\times\R^p$ we have $(x,y)\in\SP_{\{m\}}$ if and only if $m\in I_0(x,y)$
and for $(x,y)\in\SP_{\{m\}}$ we have $(x,y)\in\SP^*_{\{m\}}$ if and only if there exists
$I\in\cali(x,y)$ with $m\not\in I$.

First we prove that $\SP_{\{m\}}$ is a manifold with boundary $\SP^*_{\{m\}}$.
Indeed, since there exists a feasible point (namely $0$) for $\SQP(\bar y)$
and MFCQ holds at this point,
the constraint set $M(y)$ of $\SQP(y)$ is nonempty for any $y\in U$,
where $U\subset\R^p$ is a sufficiently small open neighborhood of $\bar y$.
Hence any problem $\SQP(y)$, $y\in U$, is solvable and its minimizer $x(y)$
is the unique stationary point of $\SQP(y)$.
It is well-known, that the mapping $y\mapsto x(y)$ is continuous,
thus $\pphi:U\to\R^n\times\R^p$, $\pphi(y):=(x(y),y)$, is a topological manifold chart
of $\SP$.

Writing elements $y$ of $\R^p$ as $\left((a_i,\alpha_i)_{i=1\dots m}\right)$,
we also use the notation $y':=\left((a_i,\alpha_i)_{i=1\dots m-1},a_m\right)$,
i.e. we have $y=(y',\alpha_m)$.
We claim the existence of an open neighborhood $V\subset\R^{n-1}$ of $\bar y'$
and a continuous mapping $\alpha_m:V\to \R$ such that
its graph $\graph(\alpha_m):=\{(y',\alpha_m(y'))~|~y'\in V\}$
parameterizes (via $\pphi$) $\SP_{\{m\}}^*$ and its epigraph
$\{(y',\alpha_m)\in V\times\R~|~\alpha_m\ge\alpha_m(y')\}$
parameterizes $\SP_{\{m\}}$.

For proving this claim, we consider the special quadratic problem
$\SQP'(y'):=\calp(f,g_1\dots g_{m-1})$
of size $(n,m-1,p)$, where $f$, $g_i$ are exactly the functions
defined in (\ref{eq:sqp}).
Note that $\SQP'(y')$ is obtained from $\SQP(y)$ by just deleting the
constraint $g_m\le 0$,
thus its constraint set contains the constraint set of $\SQP(y)$.
Hence, if the neighborhood $V\subset\R^{p-1}$ of $\bar y'$ is chosen sufficiently small,
the constraint set of $\SQP'(y')$ is also nonempty for any $y\in V$
and, as above for $\SQP$, there exists a continuous mapping
$x':V\to\R^n$ such that $x'(y')$ is the unique stationary point of $\SQP'(y')$,
$y'\in V$.
Now we are able to define $\alpha_m:V\to\R$ by setting $\alpha_m(y'):=-x'(y')^\top a_m$.
Without loss of generality we may assume $\graph(\alpha_m)\subset U$.
We have to verify the predicted properties of the mapping $\alpha_m$.

In fact, since $\SQP'(y')$ is just $\SQP(y)$ with the constraint $g_m\le 0$ deleted,
$x'(y')$ is also a stationary point for $\SQP(y)$,
provided that it is feasible for that problem,
i.e. for any $\alpha_m\le\alpha_m(y')$.
In case $\alpha_m=\alpha_m(y')$,
the constraint $g_m\le 0$ is active at $(x,y):=(x'(y'),y',\alpha_m)$.
Since $x'(y')$ is a stationary point of $\SQP'(y')$
there exists $I\in\cali(\SQP)(y)$ with $m\not\in I$,
i.e. this point $(x,y)$ belongs to $\SP_{\{m\}}^*$.
In case $\alpha_m<\alpha_m(y')$,
the constraint $g_m\le 0$ is not active at $(x,y):=(x'(y'),y',\alpha_m)$,
i.e. $(x,y)$ does not belong to $\SP_{\{m\}}$.
In case $\alpha_m>\alpha_m(y')$,
the point $x'(y')$ is not feasible for $\SQP(y)$, where $y=(y',\alpha_m)$,
since the constraint $g_m\le 0$ is violated.
Hence the stationary point $x(y)$ of $\SQP(y)$ must differ from $x'(y')$.
This implies that the constraint $g_m\le 0$ must be active at $(x(y),y)$,
since otherwise $x(y)$ would necessarily be a stationary point of $\SQP'(y')$,
thus it would be $x'(y')$.
Altogether we have seen that $\pphi:U\to\R^n\times U$ parameterizes
$\SP_{\{m\}}$ for $\alpha\ge\alpha_m(y')$ and $\SP_{\{m\}}^*$ for
$\alpha=\alpha_m(y')$,
thus $\SP_{\{m\}}$ is a manifold with boundary $\SP_{\{m\}}^*$.

For obtaining the corresponding assertion for the fibers,
it suffices to prove the existence of a local projection mapping
$\tilde\pi:\R^n\times\R^p\to\tilde Z$
such that $\tilde\pi^{-1}(0,\bar y)=\R^n\times N\times\R$
for a suitable (smooth) sub-manifold $N\subset\R^{p-1}$ containing $\bar y'$.
In fact, then $\SP\cap\tilde\pi^{-1}(0,\bar y)$
is parameterized (via $\pphi$) by
$\{(y',\alpha_m)~|~y'\in N,~\alpha_m\ge\alpha_m(y')\}$.

For any $(x,y)\in\R^n\times\R^p$ we can tell
the combinatorial code $(I_0',\cali')$ of $(x,y')$
w.r.t. $\SQP'(y')$ from the combinatorial code $(I_0,\cali)$
of $(x,y)$ w.r.t. $\SQP(y)$.
Indeed, we have $I_0'=I_0\setminus\{m\}$ and
$\cali'$ is obtained from $\cali$ by deleting all
$I\in\cali$ containing $m$.
Letting $\Pi':\R^n\times\R^p\to\R^n\times\R^{p-1}$,
$\Pi:\R^n\times\R^{p-1}\to\R^{p-1}$
denoting the canonical projections
$\Pi'(x,y):=(x,y')$, $\Pi(x,y'):=y'$,
and defining $X:=\Pi'(\tilde Z)$,
it is therefore clear,
that the combinatorial code w.r.t. $\SQP'$
is constant on $X\subset\SP(\SQP')$.
Hence we can apply Lemma \ref{lem:proj}
to $\SQP'$, $(\bar x,\bar y')\in\SP(\SQP')$,
the projection mapping $\Pi:X\to Y:=\Pi(X)$
and $U:=\R^n\times V'$,
where $V'\subset\R^{n-1}$ is an open neighborhood of $\bar y'$.
Provided that the neighborhood $V\subset V'$ of $\bar y'$ has been chosen small enough,
the assertion of Lemma \ref{lem:proj}
implies the existence of a smooth mapping $x'':V\to\R^n$
such that $x=x''(y')$ is the stationary point of $\SQP'(y')$,
i.e. $x''(y')=x'(y')$, for any $y'\in V\cap Y$.
This gives rise to a smooth mapping $\pphi'':V\to\R^n\times V$,
$\pphi''(y'):=(x''(y'),y')$.
Note that $\pphi''|_Y=\Pi|_X^{-1}$.

Moreover, for any $(x,y')\in X$, the value of $\alpha_m\in\R$
with $g_m(x,y',\alpha_m)=0$ is unique.
Since $g_m\le 0$ is active on $\tilde Z$,
we therefore have $\Pi'|_{\tilde Z}^{-1}=\pphi'|_X$,
where $\pphi'(x,y'):=(x,y',-x^\top a_m)$.

Since $\pphi':\R^n\times\R^{p-1}\to\R^n\times\R^p$ is smooth,
the mapping $\pphi'\circ\pphi:V\to\R^n\times\R^p$ is smooth
and we have $\pphi'\circ\pphi|_Y=(\Pi\circ\Pi')|_{\tilde Z}^{-1}$.
Hence $Y\cap V$ is a smooth sub-manifold of $\R^{p-1}$.

Now let $\pi_Y:\R^{p-1}\to Y$ be an arbitrary smooth projection mapping,
defined on an open neighborhood (contained in $V$) of $\bar y'$ and put
$\tilde\pi:=\pphi'\circ\pphi\circ\pi_Y\circ\Pi\circ\Pi'$.
Then we have $\tilde\pi(x,y)=(0,\bar y)$ if and only if
$\pi_Y\circ\Pi\circ\Pi'(x,y)=\bar y'$,
thus $\tilde\pi^{-1}(0,\bar y)=\R^n\times\pi_Y^{-1}(\bar y')\times\R$.
Putting $N:=\pi_Y^{-1}(\bar y')$,
we are done.\QED

\section{Final Remarks}\label{sec:final}

The following example shows that the assumption $m_==0$ of
Theorem \ref{th:manifold} cannot be omitted.
\begin{example}
\label{ex:counter}
For the problem size $(n,m)$, given by $n:=2$, $m_\le:=0$, $m_==1$,
we consider
\[
\overline{\chii\SP}=\{((b_1,\beta_1),a_{0^*})\in(\R^2\times\R)\times\R^2~|~
\beta_1=0,\ \det(b_1|a_{0^*})=0\}.
\]
The latter set, a cone over a torus, fails to have the local structure of a topological manifold
with boundary at the vertex of the cone, the point $\bar z:=((0,0),0)$.
Now we consider our standard special quadratic problem $\SQP:=\SQP_1$
from the proof of Theorem \ref{th:main} such that $\jsp(\SQP)(0,\bar y)=\bar z$.
Then $\jsp(\SQP)\transversal\{\bar z\}$ and, since $\jsp(\SQP)$ is a local
diffeomorphism, $\overline\SP$ has the same local topological structure as $\overline{\chii\SP}$.

For having this example in explicit terms, let $p:=3$ and put
\[
f(x,y):=x_1^2+x_2^2,~h_1(x,y):=x_1y_1+x_2y_2+y_3,~\bar y:=(0,0,0).
\]
\end{example}

\begin{remark}[Smoothness is not required]
The assertions of this article remain true if the assumption of smoothness
of the problem date is weakened to $C^2$.
\end{remark}

\section*{Acknowledgment}
The author is indebted to the
associate editor, Diethard Klatte, and the anonymous referees
for their precise and constructive remarks.

\end{document}